\documentclass[12pt]{article}
\usepackage{graphicx}
\usepackage{epsfig}
\usepackage{amssymb, amsmath}
\usepackage{amsthm}

\begin{document}

\title{ Spectral convergence of manifold pairs}
\author{Karsten Fissmer, Ursula Hamenst\"adt\thanks
{Research partially supported by SFB 256 and by SFB 611.}
\\
Mathematisches Institut der Universit\"at\\
Beringstra\ss{}e 1, D-53115 Bonn, Germany
}
\date{}
\maketitle

\begin{abstract}
Let $(M_i,A_i)_i$ be
pairs consisting of a complete Riemannian
manifold $M_i$ and a nonempty closed subset $A_i$. Assume that the
sequence $(M_i,A_i)_i$ converges in the Lipschitz topology to the
pair $(M,A)$. We show that there is a number $c\geq 0$ which is
determined by spectral properties of the ends of $M_i-A_i$ and
such that the intersections with $[0,c)$ of the spectra of $M_i$
converge to the intersection with $[0,c)$ of the spectrum of $M$.
This is used to construct manifolds with nontrivial essential
spectrum and arbitrarily high multiplicities for an arbitrarily
large number of eigenvalues below the essential spectrum.
\end{abstract}

\section{Introduction}

In this note we investigate the spectrum of the Laplacian acting
on square integrable functions on a complete Riemannian manifold
which is not necessarily of finite volume. Our main goal is to
understand how this spectrum varies as we vary our manifold
continuously with respect to the {\sl Lipschitz topology for
metric pairs}.

Here we mean by a {\sl metric pair} a pair $(M,A)$ which consists
of a metric space $(M,d)$ and a nonempty closed subspace $A\subset
M$. For a number $R>0$ denote by $B(A,R)$ the open
$R$-neighborhood of $A$ in $M$. The Lipschitz topology for metric
pairs is defined as follows (compare [G]).

\bigskip

{\bf Definition:} A sequence of metric pairs $(M_i,A_i)$
\emph{converges} to the metric pair $(M,A)$
\emph{in the Lipschitz topology} if
there is a sequence of numbers $R_i\to \infty$, a sequence of
numbers $\epsilon_i\to 0$ and for each $i$ a
$(1+\epsilon_i)$-bilipschitz homeomorphism
$F_i$ of
$B(A,R_i)\subset M$ onto
a neighborhood of $B(A_i,R_i)$ in $M_i$ which
maps $A$ to $A_i$.
We call the sequence $\{R_i\}_i$
\emph{convergence inducing}.

\bigskip

If the closed sets $A_i\subset M_i$ and $A\subset M$ consist of
single points then we also speak of the {\sl Lipschitz topology of
pointed metric spaces} and {\sl Lipschitz convergence of pointed
metric spaces} (see [G]).

In the sequel we only consider
metric pairs $(M,A)$ which consist of a
connected complete Riemannian manifold $M$ and a closed subset $A$
of $M$. We call such a pair $(M,A)$ a {\sl manifold pair.}

For every complete Riemannian manifold $(M,g)$, the
spectrum of the Laplacian $\Delta$ acting on square integrable
functions is a closed subset $\sigma(M)$ of the half-line
$[0,\infty)$. The set $\sigma(M)$
can be written as the union of the {\sl
essential spectrum} $\sigma_{\rm ess}(M)$ and the {\sl discrete
spectrum} $\sigma_{\rm disc}(M)$. The essential spectrum is a
closed subset of $\sigma(M)$. The discrete spectrum consists of
the eigenvalues of finite multiplicity; they are isolated points in
$\sigma(M)$.
If $M$ is closed then the essential spectrum of
$M$ is empty and $\sigma(M)$ consists of an increasing sequence
$0=\lambda_1<\lambda_2<\dots $ of nonnegative numbers converging
to $\infty$.

If $(M_i,g_i)$ are diffeomorphic closed
Riemannian manifolds which
converge as $i\to \infty$ in the Lipschitz topology to a closed
Riemannian
manifold $(M,g)$ then the spectra of $M_i$ converge to the
spectrum of $M$. However,
spectra do not always converge.
Namely, consider
a sequence $(M_i,p_i)_i$ of pointed closed
manifolds which converge in the pointed Lipschitz topology to a
complete non-compact manifold $(M,p)$ of finite volume.
Let $\nu\geq 0$ be a lower bound for the essential spectrum
of $M$ and assume that
$M$ admits exactly $k\geq 0$ eigenvalues counted with
multiplicities which are smaller than $\nu$.
In [CC1] and [CC2], Colbois and Courtois show that
the first $k$ eigenvalues of $M_i$ converge to the
first $k$ eigenvalues of $M$ if and only if there is a convergence
inducing sequence $R_i\to\infty$ such that for sufficiently
large $i$ the smallest {\sl Raleigh quotient} of $M_i-B(p_i,R_i)$
is not smaller than $\nu$. Recall that the smallest Raleigh quotient
$\mu_1(\Omega)$ of an open subset $\Omega$ of a Riemannian
manifold $(M,g)$ is defined to be the infimum of all quotients
${\cal R}(f)=\int g(df,df)/\int f^2$ over all nontrivial smooth
functions $f$ with compact support in $\Omega$.

We adapt this idea to our more general situation using the
following definition.

\bigskip

{\bf Definition:} Let $(M_i,A_i)$ be a sequence of metric pairs
converging in the Lipschitz topology to the metric pair $(M,A)$
with a convergence inducing sequence $R_i\to \infty$. A family of
open subsets $\Omega_i\subset M_i-A_i$ is called \emph{escaping} if
there is a sequence $r_i\to\infty$ such that $\Omega_i$ contains
$M_i-B(A_i,R_i-r_i)$.

\bigskip

We use here the notion of Colbois and Courtois in [CC2] even
though our definition slightly differs from theirs and our
escaping sets
do not necessarily ``escape''
in an intuitive sense.

Denote by $L^2(M)$ the Hilbert space of square integrable
functions on a Riemannian manifold $(M,g)$ and let $H^1(M)$ be the
Hilbert space of square integrable functions on $M$ with square
integrable differential. Let $(M_i,A_i)$ be a sequence of manifold
pairs converging to $(M,A)$ with convergence inducing sequence
$\{R_i\}$ and
$(1+\epsilon_i)$-bilipschitz maps $F_i:(B(A,R_i),A)\to (M_i,A_i)$.
We say that a sequence of functions $f_i\in L^2(M_i)$
{\sl converges effectively} to a function $f\in L^2(M)$ if
$\int_{M_i-B(A_i,R_i)}f_i^2\to 0$ and if moreover
$\int_{B(A,R_i)}(f_i\circ F_i -f)^2\to 0$ as $i\to \infty$. We
show.

\bigskip

{\bf Theorem A:} {\it Let $(M_i,A_i)$ be a sequence of manifold
pairs which converges in the Lipschitz topology to the manifold
pair $(M,A)$. Let $\Omega_i\subset M_i$ be an escaping family of
sets and let $\nu\leq \lim\inf_{i\to \infty}\mu_1(\Omega_i)$. Then
the sets $\sigma(M_i)\cap [0,\nu)$ converge as $i\to \infty$ in
the Hausdorff topology for closed subsets of $[0,\nu)$ to
$\sigma(M)\cap [0,\nu)$. Moreover,
every function $f\in H^1(M)$ whose spectral
measure is supported in $[0,\nu)$ is an effective limit of
functions $f_i\in H^1(M_i)$ whose spectral measures converge weakly
to the spectral measure of $f$.}

\bigskip

For closed pointed Riemannian manifolds $(M_i,p_i)$ which converge
to a complete manifold $(M,p)$ of finite volume we can combine our
Theorem A with standard compactness results for solutions of
elliptic equations to obtain convergence of eigenfunctions on
$M_i$ with small eigenvalue to eigenfunctions on $M$.

There are also sequences of closed pointed manifolds $(M_i,p_i)$
which converge in the pointed Lipschitz topology to a
noncompact manifold
$(M,p)$ and such that
up to passing
to a subsequence,
eigenfunctions on
$M_i$ with arbitrary but
bounded eigenvalues converge
locally uniformly
to an eigenfunction on $M$ (which however
may not be square integrable).
In Section 3 we look at one specific class of such examples
which can be described as follows.

Let $N$ be a closed two-sided hypersurface in a closed manifold
$M$. Then $N$ has a tubular neighborhood $U$ which is
diffeomorphic to $N\times (-1,1).$ We consider a family $g_s$ of
Riemannian metrics on $M$ which depend smoothly on $s\in (0,1]$
and which
converge as
$s\searrow 0$ uniformly on compact subsets of $M-N$ to a smooth
complete metric $g_0$ on $M-N$.
Assume that there is a smooth family $h_s$ $(s\in [0,1])$ of
smooth metrics on $N$ and a smooth function $\rho:(0,1]\times
[-1,1]\to (0,\infty)$ such that the restriction of $g_s$ to
$N\times (-1,1)$ is of the form $g_s=\frac{1}{s^2+t^2}dt^2
+\rho(s,t)h_s$.

In the following theorem we mean by an eigenfunction a solution of
the differential equation $\Delta -\lambda =0$ for some
$\lambda\in \mathbb{R}$ which need not be square integrable.

\bigskip

{\bf Theorem B:} {\it Assume that $\rho(s,t)\searrow 0$
as $(s,t)\to 0$. Let
$s_i\subset (0,1]$ be a sequence converging to $0$ and let $f_i$
be an eigenfunction on $(M,g_{s_i})$ with respect to an eigenvalue
$\lambda_i$. If the sequence $\lambda_i$ converges to some
$\lambda \geq 0$ then after passing to a subsequence and possibly
a renormalization the functions $f_i$ converge uniformly on
compact subsets of $M-N$
to a nonzero eigenfunction for $g_0$ with respect to
the eigenvalue $\lambda$.}

\bigskip

We also give a simple example for the fact that
the limit function is in general not square integrable,
even if the curvatures and the volumes of all the
metrics $g_s$
are uniformly bounded. Our theorem is due to Judge [J]
in some special cases.

In Section 4 we construct manifolds with controlled sectional
curvature and with
prescribed spectral behavior.

\bigskip

{\bf Theorem C:} {\it
For every $n\geq
2,k>0,m>0$ there is a smooth Riemannian
manifold $M$ of dimension $n$ and curvature
contained in $[-1,0]$ and with the following additional
properties.
\begin{enumerate}
\item
The essential spectrum $\sigma_{\rm ess}(M)$
of $M$ is not empty and $M$
has infinitely many
eigenvalues below $\sigma_{\rm ess}(M)$.
\item
For $2\leq j\leq k$
the multiplicity of the $j$-th eigenvalue
of the Laplacian
is at least $m$.
\end{enumerate}
In the case $n=2$ we can choose $M$ to
have constant curvature $-1$.}

\bigskip

Our construction can also be used to obtain for any $n\geq 2$
and for given $k>0,m>0$ a
{\sl compact}
$n$-dimensional manifold
with the
property that for $2\leq j\leq k$
the multiplicity of the $j$-th eigenvalue
is at least $k$. However
at least in the case of surfaces of constant curvature
such examples were known before
[BC].

\section{Proof of Theorem A}

This section is devoted to the proof of
Theorem A. We
continue to use the assumptions and notations from the
introduction. In particular, we denote by $(M,g)$ a complete
Riemannian manifold and by $A$ a nonempty closed subset of $M$.

For functions $f,h$ on $(M,g)$
denote by $(f,h)_2=\int_M fh$ their $L^2$-inner product and let
$(\nabla f,\nabla h)_2=\int_M g(\nabla f,\nabla h)$ be the
$L^2$-inner product of their gradients. Write also $\Vert
f\Vert=\sqrt{(f,f)_2}$ and $\Vert \nabla f\Vert =\sqrt{(\nabla
f,\nabla f)_2}.$ We denote by $H^1(M)$ the Hilbert space of square
integrable functions on $M$ with square integrable differential
with the inner product $(f,h)_2+(\nabla f,\nabla h)_2$. For
every nonzero function $f\in H^1(M)$ the Raleigh quotient of $f$
is defined by ${\cal R}(f)=\Vert \nabla f\Vert^2/\Vert
f\Vert^2.$

We begin our argument with a general estimate of Raleigh quotients
for suitably chosen functions on $M$.

\bigskip

{\bf Lemma 2.1:} {\it For $\epsilon >0$
there is a number
$\delta=\delta(\epsilon)>0$
with the following property. Let
$M$ be a complete Riemannian manifold,
let
$U,V$ be open subsets of $M$ with
disjoint closures and let $u,v$ be smooth functions on $M$ with
compact supports in $U,V$. If
${\cal R}(u+v)<\mu_1(V)-\epsilon$
and $\vert {\cal R}(u+v)-{\cal R}(u)\vert <\delta$ then
$\Vert
v\Vert^2<\epsilon \Vert u\Vert^2$ and
$\Vert
\nabla v\Vert^2<\epsilon(\Vert \nabla u\Vert^2
+2\epsilon \Vert u\Vert^2)$.}

{\it Proof:} Let $u,v$ be as in the lemma. Since $v$ is supported
on $V$ we have ${\cal R}(v)\geq \mu_1(V)>{\cal R}(u+v)+\epsilon$.

Write \[a=\Vert \nabla u\Vert^2,\quad b=\Vert u\Vert^2,\quad
c=\Vert\nabla v\Vert^2,\quad d=\Vert v\Vert^2.\]
Since the supports of $u$ and $v$ are disjoint
we have $\Vert u+v\Vert^2=b+d$ and $\Vert \nabla(u+v)\Vert^2=a+c$
and consequently
$\frac{c}{d }={\cal R}(v)\geq
{\cal R}(u+v)+\epsilon=
\frac{a+c}{b+d}+\epsilon$. This implies that
\[\epsilon\frac {\Vert v\Vert^2}{\Vert u\Vert^2}=
\epsilon\frac{d}{b}\leq \frac {a+c}{b+d}-  \frac{a}{b}={\cal
R}(u+v)-{\cal R}(u).\]

Thus if $\vert{\cal R}(u)-{\cal R}(u+v)\vert <
\delta(\epsilon)=\epsilon^2$
then our above inequality shows that $\Vert v\Vert^2\leq
\epsilon \Vert u\Vert^2$.

Using again that the supports of $u$ and $v$ are
disjoint we obtain from this that
\begin{equation*}
{\cal R}(u)+\epsilon^2\geq {\cal R}(u+v)\geq\frac{\Vert\nabla
u\Vert^2+\Vert \nabla v\Vert^2}{(1+\epsilon)\Vert u\Vert^2}=
\frac{{\cal R}(u)}{1+\epsilon}+ \frac{\Vert \nabla
v\Vert^2}{(1+\epsilon)\Vert u\Vert^2}
\end{equation*}
and therefore
\[\epsilon{\cal R}(u)+\epsilon^2(1+\epsilon)\geq
\frac{\Vert \nabla v\Vert^2}{\Vert u\Vert^2}\] and $\Vert\nabla
v\Vert^2\leq \epsilon \Vert \nabla u\Vert^2+
\epsilon^2(1+\epsilon)\Vert u\Vert^2.$ This shows the lemma.
{\bf q.e.d.}

\bigskip

For an open subset $\Omega$
of $M$ with smooth boundary we denote
by $H^1(\Omega)$ the closure in $H^1(M)$ of the space of smooth
functions with compact support in $\Omega$. Then $H^1(\Omega)$ is a
closed linear subspace of $H^1(M)$.

The self-adjoint extension
of the Laplacian $\Delta_\Omega$
on $\Omega$ with Dirichlet boundary conditions
is the self-adjoint operator of the quadratic form
$(f,u)\to (\nabla f,\nabla u)_2$. The domain
of $\Delta_\Omega^{1/2}$ is the
Hilbert space $H^1(\Omega)$.
We denote by $\sigma(\Omega)\subset [0,\infty)$ the spectrum
of $\Delta_\Omega$.

In the sequel we always assume without further mentioning
that the boundaries of our sets $B(A,r)$ are smooth.
This can be achieved with a small deformation
of $B(A,r)$ near its boundary. We also write $M=B(A,\infty)$.

The next lemma is the key technical result needed for the proof
of Theorem A.

\bigskip

{\bf Lemma 2.2:} {\it For $\epsilon >0,C>0,\chi \in
(0,\epsilon/2)$ there is a number $R=R(\epsilon,C,\chi)>0$ and a
number $\nu=\nu(\epsilon,C,\chi)>0$ such that the following is
satisfied. Let $M$ be a complete Riemannian manifold and let
$A\subset M$ be a closed set. Then there is a continuous linear
map $L:H^1(M)\to H^1(M)$ with the following properties.
\begin{enumerate}
\item
The range of $L$ is contained in $H^1(B(A,R))$.
\item $L$ extends continuously to $L^2(M)$, and
$(L\alpha,\beta)_2=(\alpha,L\beta)_2$ for all $\alpha,\beta\in
L^2(M)$.
\item
If $\lambda< \min\{\mu_1(M-A)-\epsilon,C\}$,
$r\in [R,\infty]$ and if the spectral
measure of $f\in H^1(B(A,r))$ is contained in
$[\lambda-\nu,\lambda+\nu]$ then $\Vert f-Lf\Vert^2 \leq \chi
\Vert f\Vert^2$ and $\Vert \nabla(Lf)-\nabla f\Vert^2 <\chi
\Vert\nabla f\Vert^2$.
\end{enumerate}
}

{\it Proof:} Let $M$ be a complete Riemannian manifold and let
$A\subset M$ be a closed set. The proof of our lemma is divided
into four steps.

\smallskip

\noindent {\sl Step 1:}

We claim that for all $C >0,\rho>0$ there is a number
$R_0=R_0(C,\rho)>0$ as follows. If $f\in H^1(M)$ is such that
$\Vert f\Vert^2=1$ and ${\cal R}(f)<C$ then there is a function
$u=u(f)$ with the following properties.
\begin{enumerate}
\item $u$ has values in $[0,1]$ and $\vert \nabla u\vert \leq 1$
pointwise.
\item $u=u_1+u_2$ where $u_1$ is supported in $B(A,R_0)$ and the
support of $u_2$ is disjoint from the support of $u_1$ and
contained in $M-A$.
\item
\[\int_{M}
\bigl((f-fu)^2+ \Vert\nabla (f-fu)\Vert^2\bigr) < \rho .\]
\end{enumerate}
For this
choose a number $k>0$
such that $k\rho/4
>1+C$. Notice that $k$ only depends on
$\rho,C$. For $m\leq k$ define $E_m=\{x\mid {\rm dist}(x,A)\in
[4m,4m+4)\}.$ Then $B(A,4k+4)-B(A,4)$ is the disjoint union of the
$k$ spherical shells $E_m$.

Let $f\in H^1(M)$ be such that $\Vert f\Vert^2=1$ and $\Vert
\nabla f\Vert^2<C$. Then $\int \Vert\nabla f\Vert^2+f^2<C+1$ and
therefore, by our choice of $k$, there is some $m\in \{1,\dots,
k\}$ such that \[\int_{E_m} \Vert\nabla f\Vert^2 + f^2<\rho/4.\]

For this number $m\leq k$, choose a smooth function $\tilde
u_1:\mathbb{R}\to [0,1]$ which is supported in $(-\infty,4m-2)$,
equals $1$ on $(-\infty,4m-4]$ and whose gradient is pointwise
bounded in norm by $1$. Similarly, let $\tilde u_2:\mathbb{R}\to
[0,1]$ be a smooth function which is supported in $(4m+2,\infty)$,
equals $1$ on $[4m+4,\infty)$ and such that the gradient of
$\tilde u_2$ is pointwise bounded in norm by $1$. Define
$u_i=\tilde u_i({\rm dist}(A,\cdot))$ and $u=u_1+u_2$. For
$R_0=4k+4$ the function $u_1$ is supported in $B(A,R_0)$. The
support of $u_2$ is contained in $M-A$ and it is disjoint from the
support of $u_1$.

The
function $1-u$ is supported in the shell $E_m$ and
it satisfies $\vert 1-u\vert \leq 1,\Vert\nabla(1-u)\Vert \leq 1$
pointwise.
Therefore we have
\begin{align*}
\int_M (f-fu)^2= \int_{E_m} f^2(1-u)^2\leq &\int_{E_m}f^2<\rho/4
\quad {\rm and}\\ \int_M\Vert\nabla (f-fu)\Vert^2= \int_{E_m}\Vert
(1-u)\nabla f + & f\nabla (1-u)\Vert^2 \\\leq
\int_{E_m}\Vert\nabla f\Vert^2 +  2(1-u)fg & (\nabla f,\nabla
(1-u))+f^2 \\<  \rho/4+ 2\int_{E_m}f\Vert\nabla f\Vert  \leq &
3\rho/4.
\end{align*}
In other words, our function $u$ has the required properties.

\smallskip
\noindent {\sl Step 2:}

We claim that for every $\epsilon >0, C>0$ there is a number
$\beta=\beta(\epsilon,C)>0$ with the following property. Let $f\in
H^1(M)$ be a normalized function with ${\cal
R}(f)<\min\{\mu_1(M-A),C\}-\epsilon$. Let $\rho\in
(0,\frac{1}{2})$ be an arbitrary number which is small enough that
$\frac{C+\rho}{1-\rho}-C-\rho< \epsilon/2$ and let $u=u_1+u_2$ be
the function constructed in step 1 for $f$ and the constants
$C,\rho >0$; then $\int (fu_1)^2\geq \beta.$

By step 1 above we have $\int f^2(1-u)^2+ \Vert \nabla
(f-fu)\Vert^2 <\rho$ and therefore since $f$ is normalized and
$u\leq 1$ we conclude that
\[ {\cal R}(f)-\rho=\frac{\Vert \nabla f\Vert^2-\rho}
{\Vert f\Vert^2} \leq R(uf)\leq \frac{\Vert \nabla
f\Vert^2+\rho}{\Vert f\Vert^2- \rho}.\] By our choice of $\rho$
and the fact that ${\cal R}(f)<C$ we conclude that $\vert {\cal
R}(uf)-{\cal R}(f)\vert <\epsilon/2$ and hence ${\cal R}(uf)<
\min\{\mu_1(M-A),C\}-\epsilon/2.$ Now
\[{\cal R}(uf)=\frac{\Vert\nabla (u_1f)\Vert^2+\Vert
\nabla(u_2f)\Vert^2}{\Vert u_1f\Vert^2+\Vert u_2 f\Vert^2}\] and
consequently since $u_2f$ is supported in $M-A$ we obtain that
\begin{align*}
\min\{\mu_1(M-A),C\}-\epsilon/2 &\\ \geq
\frac{\Vert\nabla(u_1f)\Vert^2+\mu_1(M-A)\Vert u_2f\Vert^2}{\Vert
u_1f\Vert^2+\Vert u_2f\Vert^2} &  \geq \frac{\mu_1(M-A)\Vert
u_2f\Vert^2}{\Vert u_1f\Vert^2+\Vert u_2f\Vert^2}
\end{align*} and hence
\[\Vert u_1f\Vert^2\geq \epsilon\Vert
u_2f\Vert^2/2C.\] The existence of a constant
$\beta=\beta(\epsilon,C)$ as stated above now follows from the
fact that $\Vert u_1f\Vert^2+\Vert u_2f\Vert^2\geq 1-\rho\geq
\frac{1}{2}$ by step 1.

\smallskip

\noindent{\sl Step 3:}

Let $\epsilon >0$ and let $\chi <\epsilon/2$. Let $C\geq 1$, let
$\delta=\delta(\chi/2C)<\chi/2$ be as in Lemma 2.1 and let
$\beta=\beta(\delta,C)$ be the constant from step 2 above. Notice
that $\beta$ only depends on $\epsilon,\chi,C$. Choose $\rho\in
(0,\min\{\chi/2,\delta\beta/4(3+C)\})$ small enough that $
\frac{C+\rho}{1-\rho}-C-\rho<\delta/4.$ Let $R_0=R_0(C,\rho)$ be
the constant from step 1 for $\rho$; notice that $R_0$ only
depends on $\epsilon,\chi,C$. Let $r\in [R_0+2,\infty]$ and for
simplicity write $\Omega=B(A,r)$.

We use the spectral theorem in the following form (see [D]).
There is a finite measure $\mu$ on $\sigma(\Omega)\times \mathbb{N}$
and a unitary operator $U:L^2(\Omega)\to L^2(\sigma(\Omega)\times
\mathbb{N},d\mu)$ as follows. Define $h(s,n)=s$; then
$f\in L^2(\Omega)$ lies in the domain of $\Delta_\Omega$ if and only
if $hU(f)\in L^2(\sigma(\Omega)\times \mathbb{N},d\mu)$ and
$U\Delta U^{-1}(Uf)=hU(f).$ If $f$ is contained in the domain
of $\Delta_\Omega$ then
the spectral measure of
$f$ is supported in an interval $[\lambda-\kappa,
\lambda +\kappa]$ if and
only if the function $Uf$ is supported in
$[\lambda -\kappa,\lambda +\kappa]\times \mathbb{N}$. Since
$(u,q)\to (\nabla u,\nabla q)_2$ is
the quadratic form of $\Delta_\Omega^{1/2}$ this implies
that for every $q\in H^1(\Omega)$ we have
\begin{align}
\vert (\nabla f,\nabla q)_2-\lambda (f,q)_2\vert & =
\vert \int h(Uf)(Uq)d\mu -\lambda \int (Uf)(Uq)d\mu\vert \notag\\
\leq \kappa \vert \int(Uf)(Uq)d\mu\vert & =\kappa(f,q)_2.\notag
\end{align}
Using this inequality for $u=f$ we obtain
in particular that the Raleigh quotient of $f$ is contained
in the interval $[\lambda-\kappa,\lambda +\kappa].$
Moreover, if $f$ and $q$ are contained in the domain of
$\Delta_\Omega$ and if their spectral measures are supported on
disjoint subsets of $\sigma(\Omega)$ then we have $(f,q)_2=
(\nabla f,\nabla q)_2=0.$

Let $\lambda<\min\{\mu_1(M-A),C\}-\epsilon$ and let $f\in H^1(\Omega)$
be a normalized function with spectral measure contained in
$[\lambda-\delta\sqrt{\beta}/4,\lambda +\delta\sqrt{\beta}/4]$.
Then the Raleigh quotient of $f$ is not bigger than $\lambda
+\delta/4<\min\{\mu_1(M-A),C\}-3\epsilon/8$.  Let
$u=u_1+u_2$ be the function for $f$ as in step 1 above; then as in
step 2 above we obtain that
\[ \frac{\Vert \nabla f\Vert^2-\rho}
{\Vert f\Vert^2} \leq R(uf)\leq \frac{\Vert \nabla
f\Vert^2+\rho}{\Vert f\Vert^2- \rho}\] and therefore by our choice
of $\rho$ we have $\vert R(uf)-R(f)\vert <\delta/4$ and, in
particular, $R(uf)\in [\lambda-\delta/2,\lambda+\delta/2]\subset
(0,\mu_1(M-A)-\epsilon/2]$.

On the other hand, from the properties of
the spectral measure for $f$ and the fact that $f$ is
normalized we infer that
\[\bigl \vert \int g(\nabla f,\nabla \psi)-\lambda
\int \psi f\bigr \vert \leq \delta\sqrt{\beta}\Vert\psi\Vert/4\]
for every smooth function $\psi$ on $\Omega$ with compact support.
For $\psi=u_1 f\in H^1(\Omega)$ and with the notation from step 1
above this means that
\[ \bigl |\int \Vert\nabla u_1 f\Vert^2+
\int_{E_m}g(\nabla (f(1-u_1)),\nabla (u_1f))-\lambda \int
u_1f^2\bigr | \leq \delta\sqrt{\beta} \Vert u_1f\Vert/4.\] Now
$\vert \int g(\nabla(f(1-u_1)),\nabla(u_1f))\vert\leq \int_{E_m}
\Vert \nabla f\Vert^2+f^2+2f\Vert \nabla f\Vert\leq 3\rho$ and
hence we conclude as in step 2 that
\begin{align*}
 \bigl | \int \Vert\nabla u_1f\Vert^2- \lambda \int (u_1f)^2\bigr
| & \leq 3\rho  + \lambda\int_{E_m}(1-u_1)u_1f^2+
\delta\sqrt{\beta}\Vert u_1f\Vert/4 \\ \leq
(3+\lambda)\rho+\delta\sqrt{\beta}\Vert u_1f\Vert/4 & \leq
(3+\lambda)\rho+ \delta \Vert u_1f\Vert^2/4\leq  \delta\int (u_1
f)^2/2.
\end{align*}
In particular, the Raleigh quotient ${\cal R}(u_1f)$ is contained
in $[\lambda-\delta/2,\lambda +\delta/2]$ and $\vert{\cal
R}(u_1f)-{\cal R}(uf)\vert <\delta.$

Now we can apply Lemma 2.1 to the functions $u_1f$ and $u_2f$ and
deduce that $\Vert uf-u_1f\Vert^2<\chi\Vert u_1f\Vert^2/2C$ and
$\Vert \nabla(uf-u_1f)\Vert^2<\chi \Vert \nabla u_1f\Vert^2/2C$
and therefore also $\Vert u_1f-f\Vert^2<\chi$ and $\Vert \nabla
(u_1f-f)\Vert^2<\chi$. As a consequence, we have
$\int_{M-B(A,R_0)}f^2<\chi$ and $\int_{M-B(A,R_0)}\Vert\nabla
f\Vert^2<\chi.$

\smallskip
\noindent {\sl Step 4:}

Let
$v:M\to [0,1]$ be a smooth
function with support
in $B(A,R_0+2)$ and which is constant $1$ on $B(A,R_0)$.
We may choose our function in such a way that its gradient
$\nabla v$ is pointwise bounded in norm by $1$. For a
function $f\in
H^1(M)$ define $Lf=vf$. Then $L:H^1(M)\to H^1(M)$ is clearly linear,
extends continuously to $L^2(M)$ and satisfies
$(L\alpha,\beta)_2=(\alpha,L\beta)_2$ for all
$\alpha,\beta\in L^2(M)$.
Since $\vert v\vert \leq 1 $ and $\Vert \nabla v\Vert \leq 1$
pointwise the map $L$ is continuous. More precisely, we have
$\Vert Lf-f\Vert^2 \leq \int_{M-B(A,R_0)}f^2$ and
\[\Vert \nabla
Lf-\nabla f\Vert^2 \leq \int_{B(A,R_0+2)-B(A,R_0)}f^2 +
\int_{M-B(A,R_0)}\Vert \nabla f\Vert^2.\]
This together with step 3 above shows the third
part of our lemma.

We are left with showing that the image of $H^1(M)$ under
the map $L$ is contained
in $H^1(B(A,R_0+2))$. For this observe that for every smooth
function $f$ on $M$ with compact support the function $Lf$
is smooth and compactly supported in $B(A,R_0+2)$. Since
compactly supported smooth functions are dense in $H^1(M)$ and
since $L$ is continuous,
functions with compact support in $B(A,R_0+2)$ are dense in
the range of $L$. This shows the lemma.
{\bf q.e.d.}

\bigskip

{\bf Corollary 2.3:} {\it For $\epsilon >0, C>0$ and
for $\delta <\epsilon/2$
there are
numbers $\rho=\rho(\epsilon,C,\delta) >0$ and $\kappa=\kappa(
\epsilon,C,\delta)<\delta/2$
such that for every
complete Riemannian manifold $M$ and every
closed subset $A\subset M$ the following holds.
\begin{enumerate}
\item
Let $\lambda\in [0,
\min\{\mu_1(M-A),C\}-\epsilon]\cap \sigma(M)$
and let $f\in H^1(M)$ be a function
whose spectral measure is supported in $[\lambda-\kappa,
\lambda+\kappa]$. Then there is
a function $\tilde f\in H^1(B(A,\rho))$  with spectral measure
supported in $[\lambda -\delta,
\lambda +\delta]$ and such that
$\Vert f-\tilde f\Vert^2 <\delta \Vert f\Vert^2.$
\item
Let $\lambda\in [0,
\min\{\mu_1(M-A),C\}-\epsilon]\cap \sigma(B(A,\rho))$
and let $f\in H^1(B(A,\rho))$ be a function
whose spectral measure is supported in $[\lambda-\kappa,
\lambda+\kappa]$. Then there is
a function $\tilde f\in H^1(M)$  with spectral measure
supported in $[\lambda -\delta,
\lambda +\delta]$ and such that
$\Vert f-\tilde f\Vert^2 <\delta \Vert f\Vert^2.$
\end{enumerate}
}

{\it Proof:} Let $\epsilon \in (0,1],\delta <\epsilon/2$ and let
$C\geq 1$. Define $\kappa=\delta^3/(C+1)$ and let
$\rho=R(\epsilon/2, C,\kappa^2)$ be as in Lemma 2.2.
Denote by
$L:H^1(M)\to H^1(B(A,\rho))$ the linear map from Lemma 2.2.

Let $\nu=\nu(\epsilon/2,C,\kappa^2)<\kappa/2$ be as in Lemma 2.2
and let $\lambda\in \sigma(M)\cap
[0,\min\{\mu_1(M-A),C\}-\epsilon].$ Let $f$ be a normalized
function on $M$ with spectral measure supported in
$[\lambda-\nu,\lambda+\nu]$. Then the Raleigh quotient of $f$ is
not bigger than $\lambda +\nu<\min\{\mu_1(M-A),C\}-\epsilon/2$.
Moreover, since $f$ is normalized we obtain that
\[\vert \int g(\nabla f,\nabla u)-\lambda \int uf\vert
\leq 2\nu\Vert u\Vert\] for every smooth function $u$ on $M$ with
compact support.

By construction of the operator $L$, the function $Lf$ lies in the
domain of $\Delta_{B(A,\rho)}$. Moreover by Lemma 2.2 we have
$\Vert Lf-f\Vert^2 <\kappa^2\Vert f\Vert^2$ and
$\Vert\nabla(Lf-f)\Vert^2<\kappa^2\Vert\nabla f\Vert^2$. Using the
spectral theorem for the operator $\Delta_{B(A,\rho)}$ acting on
$L^2(B(A,\rho))$ with Dirichlet boundary conditions we obtain that
the function $Lf$ admits an $L^2$-orthogonal decomposition
$Lf=\alpha +\varphi +\beta$ where the spectral measure of $\alpha$
is supported in $[0,\lambda-\delta^2]$, the spectral measure of
$\beta$ is supported in $[\lambda +\delta, \infty)$ and the
spectral measure of $\varphi$ is supported in
$[\lambda-\delta^2,\lambda +\delta]$. Since $\Vert
Lf-f\Vert^2<\kappa^2\Vert f\Vert^2$ by construction, for the first
part of our lemma it is enough to show that the square norms of
$\alpha$ and $\beta$ are bounded from above by a fixed multiple of
$\delta$.

For an
estimate of $\Vert\alpha\Vert^2$,
observe that
\[\Vert\alpha\Vert^2=\int \alpha(Lf)
=\int \alpha f+\int\alpha(Lf-f)
\leq \int \alpha f +\kappa \Vert\alpha\Vert\] since $f$ is
normalized by assumption and therefore
\begin{align}
(\lambda-\delta^2)\Vert\alpha\Vert^2 & \geq\Vert
\nabla\alpha\Vert^2 = \int g(\nabla\alpha,\nabla (Lf))  \geq  \int
g(\nabla \alpha ,\nabla f)-
\kappa\Vert\nabla\alpha\Vert\Vert\nabla f\Vert \notag\\ \geq
\lambda \int\alpha f & - \kappa (\Vert\alpha\Vert+\sqrt{C}\Vert
\nabla\alpha\Vert) \geq \lambda\Vert\alpha\Vert^2-
\kappa\Vert\alpha\Vert(\lambda +2+C).\notag
\end{align}
This shows that
$\Vert\alpha\Vert\leq 2\kappa(C+1)/
\delta^2<\delta$ by our choice of $\kappa$ and the fact that
$\lambda\leq C$.

On the other hand, the square norm of $\beta$ can be
estimated as follows. By construction and Lemma 2.2 we have
\begin{align}
(1+\kappa^2)(\lambda+\kappa)
 & \geq (1+\kappa^2)\Vert\nabla f\Vert^2\geq
 \Vert\nabla (Lf)\Vert^2 \notag\\ =
\Vert\nabla \alpha\Vert^2
&
+\Vert \nabla\varphi \Vert^2+\Vert\nabla \beta\Vert^2\geq
(\lambda -\delta^2)\Vert\varphi\Vert^2+(\lambda +\delta)
\Vert \beta\Vert^2.\notag
\end{align}
Since $\Vert\varphi\Vert^2+\Vert\beta\Vert^2=\Vert
Lf\Vert^2-\Vert\alpha\Vert^2\geq 1-\kappa^2-\Vert\alpha\Vert^2\geq
1-2\delta^2$ we obtain from this that
\[(1+\kappa^2)(\lambda +\kappa) \geq (1-2\delta^2)(\lambda -\delta^2)
+\delta\Vert\beta\Vert^2\] and hence
$\delta\Vert\beta\Vert^2\leq\kappa+\kappa^2(\lambda
+\kappa)+\delta^2+2\delta^2(\lambda-\delta^2)$ and
$\Vert\beta\Vert^2\leq \delta(3+2\lambda)$. This estimate
concludes the first part of our corollary.

To show the second part of the corollary, notice that we may
always increase $\rho$
without changing our estimates
and therefore we may assume that the first
part of our corollary is valid for $\rho$ and the constants
$\epsilon >0,C>0,\delta^2>0$. Let
$\kappa=\kappa(\epsilon,C,\delta^2)<\delta^4/8$ be the constant
from the first part of our corollary. Let $\lambda \in \sigma
(B(A,\rho))\cap [0,\min\{ \mu_1(M-A),C\}-\epsilon]$ and let $f\in
H^1(B(A,\rho))$ be a normalized function with spectral measure
supported in $[\lambda-\kappa,\lambda +\kappa]$. Then $f$ as a
function from $H^1(M)$ admits an orthogonal decomposition
$f=\alpha+\varphi +\beta$ such that the spectral measure of
$\alpha$ is supported in $[0,\lambda-2\delta^2]$, the spectral
measure of $\varphi$ is supported in $[\lambda-2\delta^2,\lambda
+\delta]$ and the spectral measure of $\beta$ is supported in
$[\lambda+\delta,\infty)$. As above it is now enough to control
the square norms of $\alpha$ and $\beta$.

For this we use our above strategy and show first
that $\Vert\beta\Vert^2\leq 3\delta+\Vert\alpha\Vert^2
(\lambda -2\delta^2)/\delta$. Namely, notice that the
functions $\alpha,\varphi,\beta$ are $L^2$-orthogonal
and also orthogonal with respect to the inner product
of $H^1(M)$. Thus the Raleigh quotient ${\cal R}(f)$ of
our function $f$
can be estimated as
\[\lambda +\kappa\geq {\cal R}(f)=\Vert\nabla \alpha\Vert^2
+\Vert \nabla\varphi \Vert^2+\Vert\nabla \beta\Vert^2\geq
(\lambda -2\delta^2)\Vert\varphi\Vert^2+(\lambda +\delta)
\Vert \beta\Vert^2.\]
Since $1-\Vert\alpha\Vert^2=\Vert\varphi\Vert^2+\Vert\beta\Vert^2$
we obtain from this that
\[\lambda +\kappa\geq (1-\Vert\alpha\Vert^2)(\lambda-2\delta^2)
+\delta\Vert\beta\Vert^2\]
and hence $\delta\Vert\beta\Vert^2\leq \kappa+2\delta^2
+(\lambda-2\delta^2)\Vert\alpha\Vert^2$
from which our above claim is immediate (recall that $\kappa\leq
\delta^4$ by assumption).

We are left with estimating $\Vert\alpha\Vert^2$. For this let
$L:H^1(M)\to H^1(B(A,\rho))$ be the operator as in Lemma 2.2.
Since the spectral measure for $f$ as a function on $B(A,\rho)$ is
contained in $[\lambda-\kappa,\lambda +\kappa]$ we deduce from
Lemma 2.2 that $\Vert Lf-f\Vert^2<\delta^2$.

The function $\alpha$ can be decomposed into a finite orthogonal
sum of functions with spectral measure supported in a subinterval
of $[0,\lambda -2\delta^2]$ of length smaller than $\kappa$. We
apply the first part of our corollary to these functions and
obtain a decomposition $L\alpha=\zeta_1+\zeta_2$ where the
spectral measure of $\zeta_1$ is supported in $[0,\lambda
-\delta^2]$ and we have $\Vert\zeta_2\Vert^2= \Vert
L\alpha-\zeta_1\Vert^2<\delta^2 \Vert\alpha\Vert^2.$ However the
spectral measure of $f$ as a function in $H^1(B(A,\rho))$ is
supported in $[\lambda-\kappa,\lambda +\kappa]$ and therefore
$\zeta_1$ is orthogonal to $f$. Thus
$(L\alpha,f)_2=(\zeta_2,f)\leq \delta\Vert\alpha\Vert$. On the
other hand,
$(L\alpha,f)_2=(\alpha,Lf)_2=(\alpha,f)_2+(\alpha,Lf-f)_2\geq
\Vert\alpha\Vert^2-\delta\Vert\alpha\Vert.$ Together with the
above this shows that $2\delta \Vert\alpha\Vert\geq
\Vert\alpha\Vert^2$ which is only possible if
$\Vert\alpha\Vert^2\leq 4\delta^2.$ Then $\Vert \beta\Vert^2\leq
\delta(4+\lambda)$ which finishes the proof of the corollary. {\bf
q.e.d.}

\bigskip

Now we are ready to show the main result of this section.

\bigskip

{\bf Proposition 2.4:} {\it Let
   $(M_i,A_i)$ be a sequence of manifold pairs
which converges
   in the Lipschitz topology to the manifold pair $(M,A)$ with
   convergence inducing sequence $R_i\to \infty$.
   Assume that there
   is
   an escaping family of sets $\Omega_i\subset M_i$
   such that
   $\smash{
           \liminf\limits_{i \to \infty} \mu_1(\Omega_i) \geq c>0}$.
Let $\sigma(M_i)\subset [0,\infty)$
be the spectrum of $M_i$ and let $\sigma(M)$ be the spectrum of
$M$. Then the following is satisfied.
\begin{enumerate}
\item
The sets $\sigma(M_i)\cap [0,c)$ converge in the
Hausdorff topology for closed subsets of $[0,c)$ to
$\sigma(M)\cap
[0,c)$.
\item
Every function $f\in H^1(M)$ with spectral measure
supported in $[0,c)$ is an effective limit of a sequence of
functions $f_i\in H^1(M_i)$ with spectral measures supported in
$[0,c)$.
\item For every $\lambda\in [0,c)\cap \sigma_{\rm disc}(M)$ and every
eigenfunction $f$ with eigenvalue $\lambda$ there is a
sequence of eigenfunctions $f_i$ on $M_i$ with respect to eigenvalues
$\lambda_i\in \sigma_{\rm disc}(M_i)\cap [0,c)$ which converge
effectively to $f$.
\end{enumerate}}

{\it Proof: } With the assumptions in the statement
of the
proposition, let $R_i\to \infty$
be a convergence inducing sequence for our convergent sequence
$(M_i,A_i)$ of manifold pairs. Choose an escaping family of sets
$\Omega_i$ with the property
that $\lim\inf_{i\to \infty}
\mu_1(\Omega_i)\geq  c>0$.
There is a sequence $r_i\to\infty$ such that
$\Omega_i\supset M-B(A_i,R_i-r_i)$.

For each $R>0$ the Laplacian acts
on the Hilbert space $H^1(B(A,R))$.
As $R\to\infty$ its spectrum $\sigma(B(A,R))$
converges
in the Hausdorff topology for closed subsets of $[0,\infty)$
to the spectrum $\sigma(M)$ of
$M$.
Since there is a
$(1+\epsilon_i)$-bilipschitz map
$F_i$ of
$B(A,R_i)$ onto a neighborhood of $B(A_i,R_i)$ in $M_i$,
this
means that as $i\to \infty$ the spectrum of the Laplacian on
$B(A_i,R_i)$ converges in the Hausdorff topology
to the spectrum of $M$.

Let $\epsilon >0$ and
for $\delta >0$ let $\rho=\rho(\epsilon/2,c,\delta)$ be
as in Corollary 2.3. If $i$ is sufficiently large then we have
$\mu_1(\Omega_i)\geq c-\epsilon/2$ and $R_i-r_i\geq \rho$.
By Corollary 2.3,
the
intersection $\sigma(M_i)\cap [0,c-\epsilon]$
is contained in the $\delta$-neighborhood of
$\sigma(B(A_i,R_i))$, and $\sigma(B(A_i,R_i))\cap[0,c-\epsilon]$
is contained in the
$\delta$-neighborhood of $\sigma(M_i)$.
Since $\epsilon >0$ and $\delta >0$ were
arbitrary we conclude that
as $i\to \infty$ (and possibly
after passing to a subsequence)
the spectrum of $M_i$ converges in the Hausdorff topology
to a closed subset $B$ of $[0,\infty)$
with the property that
$B\cap [0,c)=\sigma(M)\cap [0,c)$.
This shows the first part of our proposition.

To show the second part, let $f$ be a function on $M$ with
spectral measure supported in $[0,c-\epsilon]$. We have to show
that $f$ is an effective limit of functions on $M_i$ whose
spectral measures converge to the spectral measure of $f$. But
this follows once again from Corollary 2.3. Namely, every function
$f$ on $M$ with spectral measure contained in $[0,c-\epsilon]$ can
be approximated in $H^1(M)$ by functions supported on $B(A,R)$ for
larger and larger $R$ and with spectral measure as elements of
$H^1(B(A,R))$ supported in $[0,c-\epsilon/2]$. On the other hand,
for every $\kappa >0$, every function on $B(A,R)$ whose spectral
measure is supported in $[0,c-\epsilon/2]$ admits an orthogonal
decomposition into finitely many functions whose spectral measures
are supported on intervals of length smaller than $\kappa$. If
$\varphi\in H^1(B(A,R))$ is such a function and if $i>0$ is large
enough that $R_i>R$ then we can apply Corollary 2.3 to the
function $\varphi\circ F_i^{-1}$ on $M_i$
to obtain the required approximation.

We are left with showing
the third part of our proposition.
For this let $f$ be an eigenfunction on $M$ with eigenvalue
$\lambda \in [0,c)$. Then there is a number $\delta >0$
such that the space of functions with spectral measure
supported in $[\lambda-\delta,\lambda +\delta]$ is finite
dimensional. Our above argument immediately implies
that for sufficiently large $i$ the
dimension of the space of functions on
$M_i$ with spectral measure supported in $[\lambda-\delta/2,
\lambda +\delta/2]$ is finite as well. This completes the proof
of
our proposition.
{\bf q.e.d.}

\bigskip

For an integer $k\geq 1$ and a nonempty open subset $\Omega$ of a
Riemannian manifold the $k$-th {\sl Raleigh quotient}
$\mu_k(\Omega)$ of $\Omega$ is defined to be the infimum of all
numbers $a>0$ with the following property. For every $\epsilon
>0$ there are $k$ smooth functions
$f_1,\dots,f_k$ with compact support in $\Omega$ which are
orthonormal with respect to the $L^2$-inner product
$(f,h)_2=\int_M fh$ on $M$ and such that their gradients $\nabla
f_i$ satisfy the inequality
\[R(f_i)=\int \Vert\nabla f_i\Vert^2/\int f_i^2<a+\epsilon.\]

As an immediate corollary of Proposition 2.4 we observe.
\bigskip

{\bf Corollary 2.5:} {\it Let $(M_i,A_i)$ be a sequence of
manifold pairs converging to the manifold pair $(M,A)$ and let
$\Omega_i\subset M_i-A_i$ be a sequence of escaping sets. If
$\lim\inf _{i\to\infty}\mu_1(\Omega_i)\geq \min \sigma_{\rm
ess}(M)$ then $\mu_k(M_i)\to \mu_k(M)$ for every $k\geq 1$.}

{\it Proof:} Let $\nu_0\in [0,\infty]$ be the minimum of the
essential spectrum of $M$. If $\nu_0=\infty$ then our corollary is
immediate from Proposition 2.4, so we may assume that
$\nu_0<\infty.$ Using again Proposition 2.4 it is enough to show
that $\lim\sup_{i\to \infty}\mu_k(M_i)\leq \nu_0$ for every fixed
$k>0$. Since $\nu_0$ is contained in the essential spectrum of $M$
there is for every $k$ and every $\epsilon >0$ an orthonormal
family $f_1,\dots,f_{k}$ of functions in $L^2(M)$ with support in
a fixed compact ball $B\subset M$ and Raleigh quotients ${\cal
R}(f_i)< \nu_0+\epsilon$. For $i$ sufficiently large the set $B$
is contained in the domain of our $(1+\epsilon_i)$-bilipschitz map
$F_i$. Since $\epsilon_i\to 0$ $(i\to \infty)$ this means that for
large $i$ we can find an orthonormal family $f_1^i,\dots,f_k^i$ of
functions on $M_i$ with ${\cal R}(f_j^i)<\nu_0+2\epsilon.$ This
shows that $\lim\sup_{i\to \infty}\mu_k(M_i)\leq \nu_0.$ {\bf
q.e.d.}

\section{Development of cusps}

In this section let always $M$ be a closed
manifold
of dimension $n\geq 2$ and let $N\subset M$
be a smooth closed $2$-sided
hypersurface in $M$. Then there
is a tubular neighborhood $U$ of $N$ which is diffeomorphic to
$N\times [-1,1]$.

For $s\in [0,1]$ choose a smooth Riemannian metric $h_s$ on $N$
which depends smoothly on $s$ and let $\rho:(0,1]\times [-1,1]
\cup\{0\}\times ([-1,0)\cup (0,1])\to (0,\infty)$ be a smooth
function. Using the metrics $\rho(s,t)h_s$ on $N$ we define for
each $s> 0$ a smooth metric $g_s$ on $N\times [-1,1]$ by $g_s=
\frac{1}{t^2+s^2} dt^2 + \rho(s,t)h_s.$ As $s\searrow 0$ these
metrics converge uniformly on compact subsets of $N\times
([-1,0)\cup (0,1])$ to a a complete metric $g_0$. We assume that
the metrics $g_s$ can be extended to smooth Riemannian metrics on
$M-U$ which depend smoothly on $s\in [0,1]$. We denote these
metrics again by $g_s$, and we write $M_s$ for the manifold $M$
with the metric $g_s$ (for $s=0$ we replace $M$ by $M-N$). We
allow $M-N$ to be disconnected.

\bigskip

{\bf Lemma 3.1:} {\it The manifold pairs $(M_s,
M_s -U)$ converge as $s \to 0$
to the manifold pair $(M_0,M_0-U)$.}

{\it Proof:} By construction, the distance in $M_0$ between
the subsets $M_0-U$ and $N\times ([-\delta,0)\cup(0,\delta])$
goes to infinity as $\delta\searrow 0$.
Since by our hypothesis for $t\in [-1,0)\cup (0,1]$
the metrics $\rho(s,t)h_s$ on $N$
converge as
$s\searrow 0$ locally uniformly
in
$t\in [-1,0)\cup (0,1]$ to the metrics
$\rho(0,t)h_0$
our lemma follows.
{\bf q.e.d.}
\bigskip

Assume from now on that
the second eigenvalue of the metric $\rho(s,t)h_s$
on $N$ goes
to $\infty$ as $(s,t)\to (0,0)$. Since the metrics
$h_s$ are defined for every $s\in [0,1]$
this is equivalent to requiring
that our function
$\rho$ extends continuously to $0$ at $(0,0)$.
The following example for this situation
will be used in Section 4.

\bigskip

{\bf Example 3.2:} Let $M$ be a smooth noncompact orientable
$n$-dimensional hyperbolic manifold of finite volume. Then $M$
has a finite
number $k\geq 1$ of {\sl standard cusps}. These cusps are given by a
two-sided closed embedded hypersurface
$N\subset M$ which consists of $k$ connected components
and divides $M$ into a
manifold $\bar M$ and the cusps $E_1,\dots,
E_k$. The metric $h$ on $N$ induced
from the hyperbolic metric is flat and therefore $N$ is a finite
quotient of a collection
of $k$ tori of dimension $n-1$.
The union $\cup_{i=1}^k E_i$ of our
ends $E_1,\dots,E_k$
is diffeomorphic
to $N\times [0,\infty)$ and carries the warped product metric
$dt^2+ e^{-2t}h$.

Choose a fixed smooth convex function $\alpha:\mathbb{R}\to
(0,\infty)$ with the property that $\alpha(t)=e^{-t}$ for $t\leq
0$, $\alpha(t)= e^{-1}$ for large $t$ and such that
$\alpha^{\prime}\geq -\alpha$ and $\alpha^{\prime\prime}\leq
\alpha$. For each fixed $s>0$ define a new metric $g_s$ on
$N\times [0,\infty)$ by $g_s=dt^2+e^{-2s}\alpha(t-s)^2h$. Then the
metric $g_s$ coincides with the hyperbolic metric on $N\times [0,s]$
and extends to a complete smooth metric on all of $M$ which
coincides with the hyperbolic metric on $\bar M$.
We denote this
metric again by $g_s$.
The sectional
curvature of $g_s$
is contained in $[-1,0]$.
There is a number $\tau_0 >0$ not
depending on $s$ such that
the restriction of $g_s$ to $N\times [s+\tau_0/2,\infty)$ is
the flat product metric
$e^{2(-s-1)}h\times [0,\infty)$.

Write
$E_s=N\times (s+\tau_0,\infty)$. We can glue two copies of
$M-E_s$ along the boundary with
the natural isometry between the two boundary manifolds
$N\times \{s+\tau_0\}$
to obtain a compact
connected Riemannian manifold $M_s$. This manifold
contains two isometric copies of $\bar M$ and a totally geodesic
embedded flat hypersurface which corresponds to the boundary
components of the ends $E_s$.
If we denote by $A_s$ the union of our two copies of
$\bar M$ in $M_s$
then the manifold pairs $(M_s,A_s)$ converge in the Lipschitz
topology to the disconnected hyperbolic manifold
pair $(M_0,A_0)$ which consists of
two copies of the pair $(M,\bar M)$. We call such a
converging sequence of manifolds a {\sl standard cusp convergence}.
With respect to a suitable parametrization of the cylinders
$M_s-A_s$
in $M_s$ our family of metrics can be represented as a
$1$-parameter family of
warped product metrics
of the above form.

\bigskip

Let $\nu_s(t)$ be the volume element of the metric $\rho(s,t)h_s$
on $N$. For $s\in [0,1]$ let $W_s\subset H^1(M_s)$ be the closure
in $H^1(M_s)$ of the space of smooth functions $f$ on $M_s$ which
satisfy $\int_{N\times
\{t\}}f d\nu_s(t)=0$ for all $t\in [-1/2,1/2]$.
Denote by $\mu_s$ the volume element of the metric
$g_s$ on $M$. In the sequel we write $\int\Vert\nabla f\Vert^2d\mu_s$
to denote the integral of the square norm of the differential
with respect to the metric $g_s$.
The next lemma is similar to Lemma 2.2.

\bigskip

{\bf Lemma 3.3:} {\it For every $\epsilon >0,c>0$ there is a
number $\delta=\delta(\epsilon,c) >0$ with the following property.
Let $s\leq \delta$ and let $f\in W_s$ be a function
with
$\int_{N\times [-1/2,1/2]}
\Vert \nabla f\Vert^2d\mu_s<c
\int f^2 d\mu_s$. Then we have
\[\int_{N\times [-\delta,\delta]}f^2d\mu_s
<\epsilon \Vert f\Vert^2.\]
In particular, the Hilbert space $W_s\subset H^1(M_s)$ is compactly
embedded in $L^2(M_s)$.}

{\it Proof:}
Let $\mu_2(s,t)$ be the second Raleigh quotient of the
metric $\rho(s,t)h_s$ on $N$.
By our assumption we have $\mu_2(s,t)\to \infty$ as
$(s,t)\to (0,0)$  and therefore for every $k>0$ there is
a number $\tau=\tau(k)\in (0,1/2)$ such that
$\mu_2(s,t)>k$ for
all $s<\tau$, all $t$ with $\vert t\vert <\tau$.

Now if
$f\in W_s$
then for every $t\in [-1/2,1/2]$
the restriction of $f$ to $N\times \{t\}$
is orthogonal to the constant functions.
Moreover the measure $\mu_s$ can be represented in the
form $d\nu_s\times \alpha(s,t)dt$ for a smooth function
$\alpha\geq 1$.
Consequently
for $s<\tau$ we have
\begin{align}
\int_{N\times [-1/2,1/2]}
& \Vert \nabla f \Vert^2 d\mu_s\notag\\ \geq
\int_{-1/2}^{1/2} \bigl(\int_{N\times \{t\}}\mu_2(s,t)f^2
d\nu_s(t)\bigr) &
dt\geq k\int_{N\times [-\tau,\tau]}f^2 d\mu_s.\notag
\end{align}
If $\int_{N\times [-1/2,1/2]}
\Vert\nabla f\Vert^2d\mu_s<c\int  f^2d\mu_s$
for some $c>0$ then we deduce from
this that
$\int_{N\times [-\tau,\tau]}f^2 d\mu_s< \frac{c}{k}\int f^2d\mu_s$
which shows the first part of our lemma. Compactness of the
embedding $W_s\subset H^1(M_s) \to L^2(M_s)$ then follows from
standard compactness results. {\bf q.e.d.}

\bigskip

The following proposition generalizes an earlier result of Judge
[J] and shows our Theorem B. Its proof uses the ideas of Judge
[J], with our simple Lemma 3.3 as the main new ingredient.
In contrast to Section 2 we now mean by an eigenfunction
a solution of an equation $\Delta -\lambda=0$ for some
$\lambda\in \mathbb{R}$ which is not required to be square
integrable.

\bigskip

{\bf Proposition 3.4:} {\it Assume that as $(s,t)\to 0$ the second
eigenvalue of the metric $\rho(s,t)h_s$ tends to $\infty$. Let
$c>0$ and let $\{s_i\}_i\subset (0,1]$ be a sequence converging to
$0$. Let $f_i$ be an eigenfunction on $M_{s_i}$ with respect to an
eigenvalue $\lambda_i\leq c$. Then up to passing to a subsequence
and renormalization, the functions $f_i$ converge locally
uniformly on $M-N$ to a nontrivial eigenfunction $f$ on $M_0$ with
respect to the eigenvalue $\lambda_0=\lim_{i\to
\infty}\lambda_i.$}

{\it Proof:}
Define a linear projection $P_s:L^2(U\subset
M_s)\to L^2(U\subset M_s)$ by
\[P_sf(x,t)=\int_{N\times \{t\}}f d\nu_s(t).\]
In other words, $P_sf$ is the function which is obtained by
integration of $f$ along the manifolds $M\times \{t\}$ with
respect to the volume form of the metric $\rho(s,t)h_s$.

For $i\geq 0$ let $f_i$ be an eigenfunction on $M_{s_i}$ with
respect to the eigenvalue $\lambda_i$.
We assume that these eigenvalues are bounded from above
by some $c>0$. Let $\delta=\delta(1/2,2c)$ be as in Lemma 3.3.
Using an idea of Judge [J] we define
\begin{equation*}
\tilde f_{i}(x,t)=
\begin{cases}
f_i(x,t), &\text{if $\vert t\vert \geq \delta$;}\\
(f_i-P_{s_i}f_i)(x,t) &\text{otherwise.}
\end{cases}
\end{equation*}
To simplify our
notation we assume that the functions $\tilde f_i$ are normalized;
this only depends on the normalization of $f_i$.

Let $\alpha:(-1,1)\to [0,1]$ be a smooth function supported in
$[-3/4,3/4]$ with $\alpha(t)=1$ for $t\in [-5/8,5/8]$ and define $
u_{i}(x,t)=f_{i}(x,t)-\alpha(t) P_{s_i}f_i(x,t)$. By our
normalization assumption the $L^2$-norm of the function $ u_{i}$
is not bigger than $1$, moreover $u_i$ is contained in
$W_{s_i}$.

We claim that the $L^2$-norm of the gradient of $u_i$ is bounded
independent of $i$. To see this
recall that our metrics $g_s$ are warped product metrics on
$N\times [-1,1]$ and therefore for each fixed $s\in (0,1]$,
$t\in (-1,1)$ and every smooth function $\varphi$ on $M_s$
we have $\int_{N\times \{t\}}g_s(\nabla(\varphi-P_s\varphi),
\nabla(P_s\varphi))d\nu_s(t)=0$.

Let $\beta:(-1,1)\to [0,1]$ be a
smooth function with compact support which equals $1$ on
$[-3/4,3/4]$. Define $v_i(x,t)=\beta(t)( f_i-P_{s_i} f_i)(x,t)$;
then $\Vert v_i\Vert^2\leq 1$. Since $f_i$ is an eigenfunction
with respect to the eigenvalue $\lambda_i$, by the definition of
$v_i$ and the above we have
\begin{align}
\lambda_i\geq \lambda_i\int
v_iu_i d\mu_{s_i}= \lambda_i\int v_i f_i d\mu_{s_i}
& =\int g_{s_i}(\nabla v_i,\nabla
f_i) d\mu_{s_i}\notag \\ \geq \int_{N\times [-3/4,3/4]}\Vert \nabla
(f_i-P_{s_i}f_i)\Vert^2d\mu_{s_i}
 & \geq \int_{N\times
[-5/8,5/8]}\Vert\nabla u_i\Vert^2d\mu_{s_i}.
\end{align}

On the other hand, let $\tilde \beta:[-1,1]\to [0,1]$ be a smooth
function supported in $[-1,1/2]\cup [1/2,1]$ which is constant $1$
on $[-1,-5/8]\cup [5/8,1]$. Write $\tilde
v(x,t)=\beta(t)f_i(x,t)$. As before we deduce that $\lambda_i\geq
\lambda_i\int \tilde v_i f_id\mu_{s_i}
\geq \int_{M-N\times [-5/8,5/8]}\Vert
\nabla f_i\Vert^2d\mu_{s_i}$.
Now for $5/8\leq \vert t\vert\leq 1$ we have
\[\nabla u_i(x,t)=\nabla
f_i(x,t)-\alpha^\prime(t)P_{s_i}f_i\frac{\partial}{\partial t}-
\alpha(t)\nabla(P_{s_i}f_i)\] and therefore there is a constant
$a>0$ not depending on $i$ such that
\begin{align}
\int_{M-N\times
[-5/8,5/8]}\Vert \nabla u_i & \Vert^2d\mu_{s_i}\notag\\
\leq a\int_{M-N\times
[-5/8,5/8]}\Vert \nabla f_i\Vert^2d\mu_{s_i}
& +a\int_{N\times
[-1,-5/8]\cup[5/8,1]}f_i^2d\mu_{s_i}.\notag
\end{align}
From this and inequality (1) above we
conclude that the $L^2$-norm of the gradient of $u_i$ is bounded
independent of $i$.

We claim that after passing to a subsequence the functions $u_i$
converge in the space of locally square integrable functions on
$M_0$ to a function $u_0$ with $\Vert
u_0\Vert^2=\lim_{i\to\infty}\Vert u_i\Vert^2\leq 1$. This is
obvious if the $L^2$-norms of the functions $u_i$ converge to $0$
as $i\to \infty$, so assume that there is some $c>0$ such that
$\Vert u_i\Vert^2\geq c$ for all $i$.
Since the $L^2$-norm of the
gradient of $u_i$ is bounded independent of $i$, the Raleigh
quotients of $u_i$ are bounded independent of $i$. Lemma 3.3 then
shows that after passing to a subsequence we may assume that the
functions $ u_{i}$ converge locally in $L^2(M_0)$ to a function
$u_0$ with $\Vert u_0\Vert^2=\lim_{i\to\infty}\Vert u_i\Vert^2$.

Next we observe that after passing to a another subsequence we may
assume that the restrictions to
$N\times ([-1,-\delta]\cup [\delta,1])$ of
the functions $\tilde f_i-u_i$ converge in $L^2(M_0)$ to a
function $\chi$. Again this is obvious if the $L^2$-norm of
$\tilde  f_i-u_i$ tends to $0$ with $i$. Otherwise observe that
the function $\tilde  f_i-u_i$ can be viewed as a function on
$[-1,-\delta]\cup [\delta,1]$. Its
$L^2$-norm with respect to
a measure which is uniformly equivalent to
the standard Lebesgue measure is at most $1$. Our
above consideration implies that the $L^2$-norms of the
derivatives of $\tilde f_i-u_i$
are bounded independent of $i$. Thus we obtain
convergence from compactness of the embedding $H^1(I)\to L^2(I)$
for a compact interval $I\subset \mathbb{R}$.

Consider again inequality (1) above.
By Lemma 3.3 and our choice of $\delta$ we
either have $\Vert u_i\Vert^2<1/2$ or
$\int_{N\times [-\delta,\delta]}u_i^2d\mu_{s_i}
<\int u_i^2d\mu_{s_i}/2$ for all sufficiently large $i$.
In both cases we conclude
that $\int_{N\times [-\delta,\delta]}u_i^2d\mu_{s_i}
\leq 1/2$ for large $i$. Thus our function $u_0$
necessarily satisfies $\int_{M-N\times [-\delta,0)
\cup (0,\delta]}u_0^2d\mu_0\leq 1/2.$ Since the function
$\chi$ is supported in $M-N\times [-\delta,\delta]$ and
$\Vert \chi+u_0\Vert^2=1$ we conclude that
after passing to a subsequence
the restrictions to $M-N\times
[-\delta,\delta]$ of the functions $f_i$ converge in $L^2(M_0)$ to a
function $f_0$ with $\Vert f_0\Vert^2\in [1/2,1]$.

The function $f_i$ is a solution of an elliptic equation with
smooth coefficients. With respect to the reference metric
$g_0$ on $M-N\times (-\delta/2,\delta/2)$ the $C^2$-norms
of these coefficients are
uniformly bounded.
Since the $L^2$-norms of the restrictions to
$M_{s_i}-N\times [-\delta,\delta]$ of the functions $f_i$
are uniformly bounded as well, standard elliptic
theory implies that for every $\epsilon >0$ there is a constant
$c(\epsilon)>0$ which bounds the $C^2$-norm of the restriction of
$f_i$ to $M-N\times [-\delta-\epsilon,\delta+\epsilon]$. Thus
after passing to a subsequence the functions $ f_i$
converge locally uniformly on $M-N\times [-\delta,\delta]$
to $f_0$. This implies that $f_0$ is
a nontrivial solution
of the differential equation $\Delta_0-\lambda =0$.

Our above argument also
shows that the function $f_0$ is the
restriction to $M-N\times [-\delta,\delta]$ of an eigenfunction on $M_0$
which is a locally uniform limit of a subsequence of our functions
$f_i$. Namely, for $k>-\log\delta +\log 2 $ define
\begin{equation*}
\tilde f_{i,k}(x,t)=
\begin{cases}
f_{i}(x,t), &\text{if $\vert t\vert \geq 2^{-k}$;}\\
(f_{i}-P_{s_i}f_i)(x,t) &\text{otherwise}
\end{cases}
\end{equation*}
and write $a_{i,k}=1/\Vert\tilde f_{i,k}\Vert$. For each fixed $i$
the sequence $\{a_{i,k}\}_k$ is monoto\-nously decreasing. As
before we conclude that after passing to a subsequence the
restrictions of $a_{i,k} f_{i,k}$ to $M-N\times [2^{-k},2^{-k}]$
converge locally uniformly to a solution $\bar f_{0,k}\not\equiv 0
$ of the equation $\Delta_0-\lambda =0$. Its restriction to
$M-N\times [-2^{-k},2^{-k}]$ necessarily coincides with a
nonnegative
multiple of $\bar f_{0,k-1}$.
Since no nontrivial solution of the equation $\Delta_0 -\lambda=0$
can vanish on a nontrivial open set the restriction of our function
$M-N\times [-2^{-k},2^{-k}]$ is in fact a positive multiple
of $\bar f_{0,k-1}$.
With a standard diagonal sequence
argument we conclude from this that after passing to a subsequence
our eigenfunctions $f_i$ converge locally uniformly to an
eigenfunction $f_0$ on $M_0$. {\bf q.e.d.}

\bigskip

The following example shows that the limit function obtained in
Proposition 3.4 is in general not square integrable, even if the
curvature of all our manifolds as well as their volumes are
uniformly bounded.

\bigskip

{\bf Example 3.6:} Consider a closed hyperbolic surface
$S$ of genus $2$ which consists of two bordered tori
$T_1,T_2$ glued at the boundary.
Choose a simple closed geodesic
$\gamma$ on $T_1$ which cuts $T_1$ into
a pair of pants. We denote by $g_s$ the hyperbolic metric on $S$
which we obtain
by leaving the torus $T_2$ and the twist
parameters for the glueings fixed and replacing the torus
$T_1$ by a torus for which the length of the geodesic $\gamma$
equals $s$. For a fixed point $q\in T_2$ the
pointed surfaces $((S,g_s),q)$
degenerate as $s\searrow 0$
to a twice punctured hyperbolic torus
$(S_0,g_0)$ with two finite volume cusps.
The essential spectrum of $S_0$ is bounded from below
by $1/4$ and the
second Raleigh quotient $\mu_2(S_0)$ of $S_0$ is positive.
The metrics $g_s$ are warped product metrics in a
tubular neighborhood
of $\gamma$.

Choose a number $k>0$ such that there is a smooth nontrivial
compactly supported function $f$ on the interval $(0,k)$ which
satisfies $\int_0^k f=0$ and $\int_0^k (f^\prime)^2
<\mu_2(S_0)\int_0^kf^2 /2.$ For $a>0$ and $\tau\in [0,k]$ denote
by $C_{a,\tau}$ the cylinder $S^1\times [0,\tau]$ with the metric
$a^2 ds^2+dt^2$ (where $ds^2$ is the length element of total
length $1$ on $S^1$). For every $a>0$ the function $f$ can be
viewed as a function on the cylinder $C_{a,k}$ which only depends
on the second coordinate. We have $\int_{C_{a,k}}f=0$ and
$\int_{C_{a,k}}\Vert\nabla f\Vert^2
<\mu_2(S_0)\int_{C_{a,k}}f^2/2$ for all $a>0$.

For $s\in (0,1]$ and $\tau\in [0,k]$ we replace the metric $g_s$
near $\gamma$ by a metric $\tilde g_{s,\tau}$ which is obtained
from $g_s$ by cutting $S$ open along $\gamma$ and inserting the
cylinder $C_{s,\tau}$. We slightly modify the resulting metric
near the boundary of $C_{s,\tau}$ in such a way that we obtain a
smooth metric $\tilde g_{s,\tau}$ depending smoothly on $s,\tau$
and such that $\tilde g_{s,0}=g_s$. We may assume that there is an
tubular neighborhood $Z\sim S^1\times [-1,1]$ about $\gamma$ in
$S$ such that the restrictions of the metrics $\tilde g_{s,\tau}$
to $Z$ are warped product metrics which degenerate as $s\searrow
0$ to the metric $g_0$. The metrics can be constructed in such a
way that their curvature is contained in $[-1,0]$ and that their
volumes are uniformly bounded.

For fixed $s>0$, the second Raleigh quotient of $\tilde
g_{s,\tau}$ depends continuously on $\tau\in [0,k]$. For $\tau=k$
this Raleigh quotient is not bigger than $\mu_2(S_0)/2$. Moreover,
if $s$ is sufficiently small then the second Raleigh quotient of
$\tilde g_{s,0}$ equals at least $3\mu_2(S_0)/4$ [CC2]. Thus there
is some $\tau(s)\in [0,k]$ such that this Raleigh quotient equals
exactly $\mu_2(S_0)/2$. We may assume that $\tau(s)$ depends
continuously on $s$. Define $h_s=\tilde g_{s,\tau(s)}$. Then there
is an eigenfunction $\phi_s$ on $(S,h_s)$ with respect to the
eigenvalue $\mu_2(S_0)/2$. Moreover the metrics $h_s$ satisfy the
assumptions in Proposition 3.4.

By Proposition 3.4, after renormalization and passing to a
subsequence we may assume that the eigenfunctions $\phi_s$
converge uniformly on compact sets to an eigenfunction $\phi$ on
$S_0$ with respect to the eigenvalue $\mu_2(S_0)/2$. But then
$\phi$ can not be square integrable.

\bigskip
{\bf Remark:} The considerations in Example 3.6 can also
be used to construct
for every noncompact hyperbolic surface $S$ of finite
volume and every $\lambda\in (0,\mu_2(S))$ an eigenfunction
$\phi$ on $S$ 
with respect to the eigenvalue $\lambda$. This function $\phi$ is not
square integrable.

\section{Manifolds with prescribed
spectral properties}

In this section we apply the results from Section 2 to construct
complete Riemannian manifolds of an arbitrary dimension $n\geq 2$,
with curvature in the interval $[-1,0]$, infinite volume, nonempty
essential spectrum, infinitely many eigenvalues below the
essential spectrum and with arbitrarily large multiplicities  of
an arbitrary finite number of eigenvalues. In the case $n=2$ we
can choose our manifolds to be of constant curvature $-1$.

Our manifolds will be constructed from building blocks which
consist of complete manifolds of curvature contained in $[-1,0]$
with a fixed even number $2k\geq 2$ of standard constant curvature
cusps as in Example 3.2. We describe these building
blocks in the next lemma which is a modified version of Example 4.1
of [BCD].

\bigskip

{\bf Lemma 4.1:} {\it For every $n\geq 2$,$k\geq 1$ there is a
complete $n$-dimensional Riemannian manifold $X$ of infinite
volume with the following properties.
\begin{enumerate}
\item The curvature of $X$ is contained in $[-1,0]$.
\item $X$ has $2k$
standard cusps of curvature $-1$ which are mutually isometric.
\item The essential spectrum $\sigma_{\rm ess}(X)$ of $M$ is not empty,
and there are infinitely many different eigenvalues below
$\sigma_{\rm ess}(X).$
\end{enumerate}
}

{\it Proof:} Let $\Gamma\subset SO(n,1)$ be a non-uniform lattice.
Then $V_0={\bf H}^n/\Gamma$ is a hyperbolic manifold of finite
volume with at least one end $C$. This end is a standard cusp.

The group $\Gamma$ is residually finite and therefore there is a
finite covering $V_1$ of $V_0$ such that the cusp $C$ has at least
$2k$ preimages in $V_1$. We choose $2k$ of these preimages and
label them by $C_1,\dots,C_{2k}$. If $C_1,\dots,C_{2k}$ are the
only cusps of $V_1$ then we define $N_0=V_1$.

Otherwise, i.e. if $V_1$ has additional cusps, then we modify the
metric on each of these additional cusps as in Example 3.2. These
cusps then become flat cylinders which we cut along a totally
geodesic hypersurface. We obtain a manifold $V_1^\prime$ with $2k$
cusps and a finite number of totally geodesic boundary components.
Choose a second copy $V_1^{\prime\prime}$ of $V_1^\prime$ and glue
$V_1^{\prime\prime}$ to $V_1^\prime$ along corresponding boundary
components. The resulting manifold $V_2$ is connected, its
curvature is contained in the interval $[-1,0]$ and it has
precisely $4k$ mutually isometric ends which are the cusps
$C_1,\dots,C_{2k}$ of $V_1$ and the corresponding cusps
$C_1^{\prime\prime},\dots, C_{2k}^{\prime\prime}$ of
$V_1^{\prime\prime}$. Since the cusps $C_i^{\prime\prime}$ are
mutually isometric we can replace them as before by isometric
cylindrical ends which we cut and glue in pairs to $k$ compact
handles. We obtain a manifold $N_0$ with precisely $2k$ ends. It
carries a complete Riemannian metric of finite volume with
curvature in the interval $[-1,0]$ in such a way that each of its
ends is isometric to a fixed standard cusps of constant curvature
$-1$.

Let $F_k$ be the free group with $k$ generators $\gamma_1,\dots,
\gamma_k$. We label each of the $2k$ ends of $N_0$ by one of
elements $\gamma_1,\dots,\gamma_k,\gamma_1^{-1},\dots,
\gamma_{k}^{-1}$ of $F_k$. For $s>0$ replace the standard cusps of
$N_0$ by a compact end with boundary equipped with the metric
$g_s$ from Example 3.2. The resulting manifold $N_s$ has $2k$
totally geodesic boundary components which we label as before by
the generators of $F_k$. Choose one copy of $N_s$ for every
element of $F_k$ and label it by this group element. Glue the
boundary component with label $\gamma_i$ of the copy of $N_s$ with
label $\psi\in F_k$ to the boundary component with label
$\gamma_i^{-1}$ of the copy of $N_s$ with label $\gamma_i\psi$
with the obvious isometry. We obtain a smooth manifold $M$ with a
complete Riemannian metric $g_s$ of curvature contained in
$[-1,0]$ and which depends smoothly on $s\in (0,1]$. The free
group $F_k$ acts freely and properly discontinuously on $M$ by
right translations on the labels of our basic building blocks. The
metrics $g_s$ are invariant under this action. The quotient
$(M,g_s)/F_k$ is compact and can be obtained from $N_s$ by glueing
pairwise the boundary components.

Since $(M,g_s)$
admits a free and properly discontinuous isometric action of $F_k$
with a compact quotient,
the discrete spectrum of $(M,g_s)$ vanishes and its essential
spectrum is bounded away from $0$. The bottom $\nu_s$ of this
essential spectrum depends continuously on $s$ and goes to $0$ as
$s\to 0$.

Following [BCD] we fix a number $\tau>0$ such that
$\nu_{\tau}<(n-1)^2/4$  and a sequence $\tau_i\subset
(0,\tau)$ such that $\tau_i<\tau_j$ and
$\nu_{\tau_i}<\nu_{\tau_j}$ for $i<j$. We use this
sequence to construct inductively our building block.

There is a natural word norm on the group
$F_k$ defined by our choice of generators.
For $m\geq 1$ we denote by $B(m)$ the connected
submanifold of $M$ which
consists of precisely those copies of our
manifold $N_s$ which are labeled by elements of $F_k$ of
word norm
at most $m$.
Then $B(m)$ is a smooth submanifold of $M$ with boundary. Each
of its boundary components
is totally geodesic with respect to $g_s$.
The set $B(0)$ is just the copy
of $N_s$ which corresponds to the unit element in $F_k$.

In our first step we determine a number $m_1>0$ such that there is
a function $\psi_1$ on $(M,g_{\tau_1})$ which is supported in
$B(m_1-1)-B(1)$ and with Raleigh quotient ${\cal R}(\psi_1)<
\nu_{\tau_2}$.
Modify the metric of $B(m_1)$ near
the boundary so that the new metric coincides
with $g_{\tau_1}$ on $B(m_1-1)$ and
with $g_{\tau_2}$ near the boundary.
Glue the resulting
manifold along its boundary to $(M,g_{\tau_2})-B(m_1)$. We obtain
a new manifold $\tilde M_1$. Since the essential spectrum of any
Riemannian manifold does not change under a compactly supported
change of the metric, the bottom of
the essential spectrum of $\tilde M_1$ equals $\nu_{\tau_2}$.
But the Raleigh
quotient of the function $\psi_1$ on
$B(m_1-1)\subset \tilde M_1$
is smaller than $\nu_{\tau_2}$
and hence
$\tilde M_1$ has an eigenvalue below its essential
spectrum.

We can now iterate this construction.
In our $i$-th step we begin
with a metric $\tilde g_i$ on $M$ which coincides with the metric
$g_{\tau_i}$ on $M-B(m_i)$ for some $m_i>0$ and such that there are
$i$ functions on $(M,\tilde g_i)$ with pairwise disjoint support
contained in $B(m_i-1)-B(1)$ and with Raleigh quotients smaller than
$\nu_{\tau_{i}}$.
There are at least
$i$ distinct eigenvalues below the essential spectrum. Choose a
function $\psi_{i+1}$ supported on $B(m_{i+1}-1)-B(m_i)$ for some
$m_{i+1}>m_i$ with Raleigh quotient smaller than
$\nu_{\tau_{i+1}}$. As before we change the metric outside
$B(m_{i+1})$ to $g_{\tau_{i+1}}$ to obtain a new metric $\tilde
g_{i+1}$ which admits at least $i+1$ distinct eigenvalues below the
essential spectrum.

We can repeat this construction infinitely often to obtain a
complete manifold $X_0$ with infinitely many eigenvalues below the
essential spectrum. The lower bound $\nu_0$ of the essential spectrum
of $X_0$ is strictly smaller than $(n-1)^2/4$.

Remove $B(0)$ from $X_0$ and replace each of the
$2k$ boundary components of the resulting
manifold by a standard cusp. We claim that the complete
Riemannian manifold $X$ which we obtain in this way
has the properties stated in our lemma. To see this recall
that the bottom of the essential spectrum of a standard
hyperbolic cusp
equals $(n-1)^2/4>\nu_{0}$. Since the essential spectrum
of a complete Riemannian manifold equals the essential  spectrum
of its ends, the bottom of the essential
spectrum is $\nu_0$.
The functions on $X_0$ which we constructed
above are supported outside $B(1)$
and hence can be viewed as functions on $X$. This implies that
there are infinitely many
distinct eigenvalues below the bottom
of the essential spectrum on $X$.

Using pairs of pants as in [BCD] it is clear that
for $n=2$ we can choose our manifold to have
constant curvature $-1$.
{\bf q.e.d.}

\bigskip

Consider now an arbitrary {\sl finite} group $\Gamma$. We call
$\Gamma$ {\sl admissible} if $\Gamma$ can be generated by elements
of order at least $3$. A set of generators $\gamma_1,\dots
,\gamma_{2k}$ of $\Gamma$ is called {\sl admissible} if it
consists of elements of order at least $3$, contains with each
element also its inverse and is minimal with this property.

Let $G$ be the Cayley graph for $\Gamma$ with respect to our
generators. Then $G$ is a finite connected graph whose vertices
correspond to the elements of $\Gamma$. By our choice of
generators the graph $G$ is {\sl simple} (i.e. no multiple edges
and no loops) and $2k$-regular [dH]. Two vertices $a,b\in \Gamma$
of $G$ are connected by an edge if and only if there is some $i$
such that $b=\gamma_i a$. Right multiplication in $\Gamma$ induces
an action of $\Gamma$ as a group of automorphisms of $G$ which is
transitive on the vertices.

Assume that $\gamma_{2i}=\gamma_{2i-1}^{-1}$ for $1\leq i\leq k$.
Let $X$ be a manifold as in Lemma 4.1 with $2k$ standard cusps.
We label each of these cusps by one of our generators $\gamma_i$
of $\Gamma$. For a
$k$-tuple $a=(s_1,\dots,s_k)\in [0,1]^k$
we construct a complete Riemannian manifold
$M(a)$ as follows: Choose $\vert \Gamma
\vert$ copies of $X$ and label
each of these copies with a different element of $\Gamma$. For
$1\leq i\leq k$
replace the standard cusps
of $X$ which are labeled by $\gamma_{2i-1},\gamma_{2i}$ by a
compact end with boundary equipped with the metric $g_{s_i}$. The
boundary components of the resulting manifold correspond to our
generators $\gamma_1,\dots,\gamma_{2k}$. Glue the boundary
component $\gamma_i$ of the copy of $X$ with label
$\psi\in \Gamma$ to
the boundary component $\gamma_i^{-1}$ of the copy of $X$ with
label $\gamma_i \psi$ by the obvious isometry
as before.
We obtain a complete Riemannian manifold
$M(s_1,\dots,s_k)$ which consists of $\vert \Gamma\vert$ copies of
$X$ glued at their boundaries. It contains a
distinguished collection of totally geodesic embedded
hypersurfaces, and its curvature is contained in $[-1,0]$. The
essential spectrum of $M(s_1,\dots,s_k)$ is bounded away from $0$
and there are infinitely many
eigenvalues below the essential spectrum.
The manifold $M(0)$
consists of $\vert \Gamma\vert $ copies of $X$.
We call the manifold
$M(s_1,\dots,s_k)$ a {\sl $\Gamma$-graph manifold}, and its metric
a {\sl $(s_1,\dots,s_k)$-graph metric}.

\bigskip

{\bf Lemma 4.2:} {\it For each $a\in [0,1]^k$ the group $\Gamma$
acts freely and isometrically on $M(a)$. For every fixed $q\geq 1$
the assignment $a\in [0,1]^k\to \mu_q(M(a))$ is continuous.}

{\it Proof:} Every element of $\Gamma$ acts on
the Cayley graph $G$ by an
automorphism which permutes the edges with a given label.
For each $a\in [0,1]^k$
this automorphism induces an isometry of our manifold $M(a)$
which permutes our copies of $X$ and
preserves each of the $k$
collections of hypersurfaces corresponding to one of our
generators or its inverse.
Since the action of $\Gamma$ on $G$ is
free the same is true for the action of $\Gamma$ on $M(a)$.
Continuity of the assignment $a\in [0,1]^k\to \mu_q(M(a))$
is immediate from Corollary 2.5.
{\bf
q.e.d.}

\bigskip

Let again $\Gamma$ be an admissible finite group with an
admissible set $\gamma_1,\dots,\gamma_{2k}$ of generators and
corresponding Cayley graph $G$. By
definition, this set of generators is minimal with the property
that it consists of elements of order at least $3$ and contains
with each element its inverse. Thus if we fix some $i\leq 2k$ and
if we delete all the edges in $G$
which are either labeled by
$\gamma_i$ or by $\gamma_i^{-1}$ then the resulting graph is
disconnected.

Recall that for every complete Riemannian manifold which has
eigenvalues below the essential spectrum the multiplicity of the
smallest eigenvalue is one. Following the basic idea of [BC]
we use isometric actions of finite groups
to construct complete manifolds of bounded curvature with
infinitely many eigenvalues below the essential spectrum and such
that the multiplicity of the second eigenvalue is bigger than $1$.

\bigskip

{\bf Lemma 4.3:} {\it Let $\Gamma$ be an admissible group with an
admissible set of generators $\gamma_1,\dots,\gamma_{2k}$. Let
$m\geq 2$ be the minimal dimension of a nontrivial irreducible
orthogonal representation of $\Gamma$. Then for every $a\in
(0,1]^k$ which is sufficiently close to $0$, the multiplicity of
the second eigenvalue of the $\Gamma$-graph manifold $M(a)$
is at least $m$.}

{\it Proof:} By our assumption, for each $a\in [0,1]^k$ the
quotient $M(a)/\Gamma$ is a complete manifold, and the
projection $M(a)\to M(a)/\Gamma$ is a smooth $q$-sheeted covering
where $q$ is the cardinality of $G$.

The pair
$\gamma_{2i-1},\gamma_{2i}$
of generators of $\Gamma$
defines a $\Gamma$-orbit of edges in the
Cayley graph $G$. If $s_i>0$ then this orbit of edges
corresponds to a $\Gamma$-orbit of totally geodesic
embedded closed hypersurfaces in $M(a)$ which projects to a
closed totally geodesic embedded non-separating hypersurface
in $M(a)/\Gamma$. As $s_i\searrow 0$ this hypersurface
in $M(a)/\Gamma$ degenerates to a pair of cusps.

For each $a\in (0,1]^k$ the manifold $M(a)$ is connected.
The bottom of the spectrum of $M(a)$ is not
contained in the essential spectrum and therefore it is an
eigenvalue of multiplicity $1$.
The isometric action of $\Gamma$ on
$M(a)$ induces an orthogonal representation of $\Gamma$
on
the corresponding eigenspace.
Since the dimension of this eigenspace is $1$,
this representation
is trivial and every
eigenfunction
with respect to this eigenvalue is
$\Gamma$-invariant and projects
to an eigenfunction on $M(a)/\Gamma$. In particular, the smallest
eigenvalue of $M(a)$ coincides with the smallest eigenvalue of
$M(a)/\Gamma$.

Now let $a=(0,s_2,\dots,s_k)\in [0,1]^k$
where $s_i>0$ for $i\geq 2$.
By minimality of our set of generators for $\Gamma$,
$M(a)$ consists of
at least two isometric components which are permuted by the action
of the group $\Gamma$. Thus the multiplicity of the
smallest eigenvalue of $M(a)$ (which equals the number of connected
components of $M(a)$) is at least $2$.

By Theorem A from the introduction,
as $s\searrow 0$ the
eigenvalues of $Q(s)=M(s,s_2,\dots,s_k)$ converge to the
eigenvalues of $M(a)$. The multiplicity of the first
eigenvalue of $Q(s)$ is $1$ and hence
for sufficiently small $s$ the second
eigenvalue of $Q(s)$ is strictly smaller than the second
eigenvalue of $M(a)/\Gamma$. Then
an eigenfunction for
this eigenvalue can not be $\Gamma$-invariant. This means that the
natural orthogonal representation of $\Gamma$ on the eigenspace
of $Q(s)$
with respect to the second eigenvalue
does not contain a trivial
component and the dimension of this
eigenspace equals at least the minimal dimension of a
nontrivial irreducible orthogonal
representation of $\Gamma$. This finishes the
proof of our lemma. {\bf q.e.d.}

\bigskip

We can now iterate this construction as follows. Assume that
$\Gamma$ is a finite group which contains a nested sequence
$\Gamma\triangleright H_1\triangleright\dots \triangleright H_m$
of admissible subgroups $H_i$.

Define a set of
generators $\gamma_1,\dots,\gamma_{2k}$ of $\Gamma$ to be {\sl$(
\Gamma,H_1,\dots,H_m)$-admissible} if the following is satisfied.
\begin{enumerate}
\item For every $i\leq m$ there is some $j(i)<k$ such that
$\gamma_1,\dots,\gamma_{2j(i)}$ is an admissible set of generators
for $H_i$.
\item For each $i$ the subgroup of $\Gamma$ which is generated by
those of our generators which are not contained in $H_i$
intersects $H_i$ only in the unit element.
\end{enumerate}

We call $(\Gamma,H_1,\dots,H_m)$ an {\sl admissible sequence of
groups} if
it admits a $(\Gamma,H_1,\dots,H_m)$-admissible set of generators
and if moreover
for every $i\geq 1$
the group $H_{i+1}$ is a proper normal subgroup of $H_i$.
We do not require that $H_1$ is a normal subgroup of $\Gamma$.

Now let $(\Gamma,H_1,\dots,H_m)$ be
an admissible sequence of groups.
For a given choice of a
basic manifold $X$ with $2k$ standard cusps
as in Lemma 4.1 we
constructed above from the Cayley graph of an
$(\Gamma,H_1,\dots,H_m)$-admissible set
$\gamma_1,\dots,\gamma_{2k}$ of generators
 a connected smooth manifold
$M$ which admits a natural free action of $\Gamma$ by
diffeomorphisms and a natural family of $\Gamma$-invariant
metrics.
We call our manifold $M$
a {\sl $ (\Gamma,H_1,\dots,H_m)$-graph manifold}.

\bigskip

{\bf Corollary 4.5:} {\it Let $(\Gamma,H_1,\dots,H_m)$
be an admissible sequence of groups.
Let $q\geq 2$
be the minimal dimension of
an irreducible orthogonal representation of $\Gamma$ whose
restriction to $H_1$ is non-trivial. Then there is
a family of
$(\Gamma,H_1,\dots,H_m)$-graph manifolds
for which the multiplicity of the $j$-th
eigenvalue for $j=2,\dots, m$ is at least $q$.}

{\it Proof:} Let $M$ be a $(\Gamma,H_1,\dots,H_m)$-graph
manifold.
The group $\Gamma$ acts on $M$
freely as a group of diffeomorphisms and $M\to M/\Gamma=M_0$ is a
$\vert \Gamma\vert$-sheeted covering.
Every complete metric on $M_0$ lifts to a
$\Gamma$-invariant complete metric on $M$.

Let $\gamma_1,\dots,\gamma_{2k}$ be a
$(\Gamma,H_1,\dots,H_m)$-admissible set of generators for
$\Gamma$. Let $\ell\leq k$ be such that the set
$\gamma_1,\dots,\gamma_{2\ell}$ generates $H_1$. Denote by $E$ the
subgroup of $\Gamma$ generated by
$\gamma_{2\ell+1},\dots,\gamma_{2k}$. Then $E$ is an admissible
finite group which intersects $H_1$ trivially. The Cayley graph
$G^\prime$ of $E$ with respect to the generators $\gamma_{2\ell
+1},\dots,\gamma_{2k}$ is a connected subgraph of $G$. If we
remove from $G$ all the edges corresponding to the generators
$\gamma_1,\dots,\gamma_{2\ell}$ then the resulting graph consists
of $\vert H_1\vert$ disjoint copies of $G^\prime$.

Fix some $s_0>0$ and
for $(s_1,\dots,s_{m})\in [0,1]^{m}$ let
$Q(s_1,\dots,s_m)$ be the $\Gamma$-graph manifold
$M(a)$ where $a=(a_1,\dots,a_k)$ is the $k$-tuple
defined as follows: For each $j$ let $i\leq j$ be such
that the generator $\gamma_j$ is contained in the group
$H_i$ but not in the group $H_{i+1}$ (where we put $H_0=
\Gamma$) and define $a_j=s_i$.
The group $H_1$ acts on $Q(s_1,\dots,s_m)$ as a group of
isometries. The manifold
$Q(0,\dots,0)/H_1$ is connected.

Now apply the considerations in the proof of Lemma 4.4 to the
graph manifolds $Q(s,0,\dots,0)$. Since the subgroup of $\Gamma$
generated by those elements of $\gamma_1,\dots,´\gamma_{2k}$ which
are not contained in $H_2$ intersects $H_2$ trivially, the
covering $Q(0,\dots,0)/H_2$ of $Q(0,\dots,0)/H_1$ with deck group
$H_1/H_2$ is disconnected and for $s>0$ the covering
$Q(s,0,\dots,0)/H_2$ of $Q(s,0,\dots,0)/H_1$ with deck group
$H_1/H_2$ is connected. By the considerations in the proof of
Lemma 4.4, for sufficiently small $s$ the quotient
$Q(s,0,\dots,0)/H_2$ is a complete connected manifold with
$H_1/H_2$-invariant metric for which the second eigenvalue is
strictly smaller than the second eigenvalue of
$Q(s,0,\dots,0)/H_1$. Proceeding inductively we obtain in $m-1$
steps a $(\Gamma,H_1,\dots,H_m)$-graph manifold of the form
$Q(s_1,\dots,s_m)$ such that the $j$-th eigenvalue for $2\leq
j\leq m$ is strictly smaller than the second eigenvalue of the
quotient $Q(s_1,\dots,s_m)/H_1$. This then implies that the
representation of $\Gamma$ on each of the corresponding
eigenspaces does not factor to a representation of $\Gamma/H_1$.
This shows the corollary. {\bf q.e.d.}

\bigskip

It remains to find admissible finite groups $\Gamma$ with
arbitrarily long nested sequences of admissible subgroups $H_i$
and for which the smallest dimension of an irreducible
representation which is nontrivial on $H_1$
is arbitrarily large. This is satisfied for the groups
which were already considered by Burger and Colbois [BC].

Namely, let $p\geq 3$ be an odd prime and for some $n\geq 1$ let
$\mathbb{F}_q$ be the field with $q=p^n$ elements and
multiplicative group $\mathbb{F}_q^*=\mathbb{F}_q -\{0\}$. For a
divisor $r$ of $n$ write $m=(p^n-1)/(p^r-1)$ and define
\begin{equation*}
G_{q,m}=\left\{ \left(
\begin{matrix}
\alpha^m &\beta\\
0 & 1
\end{matrix}
\right) \mid \alpha\in \mathbb{F}_q^*,\beta
\in \mathbb{F}_p\right\}.
\end{equation*}
Then $G_{q,m}$ is the semi-direct product of $\mathbb{F}_q$ with
the cyclic group $\mathbb{A}_{q,m}=\{a^m\mid a\in
\mathbb{F}_q^*\}$ of order $p^r-1$ which acts on $\mathbb{F}_q$ by
multiplication. Its commutator subgroup $H_1$ is the cyclic group
\begin{equation*}
H_{1}=\left\{ \left(
\begin{matrix}
1 &\beta\\
0 & 1
\end{matrix}
\right) \mid \beta \in \mathbb{F}_q\right\}
\end{equation*}
of order $q$ which can naturally be identified with
the additive group $\mathbb{F}_q$.

Let $\xi$ be a generator of the cyclic group $\mathbb{F}_q^*$.
Then
the
dimension of $\mathbb{F}_q$ as a vector space over the field
$\mathbb{F}_q[\xi^m]$ equals $n/r$. Choose a basis
$g_1,\dots,g_{n/r}\subset \mathbb{F}_q$ for this vector space.
For each $i$ the element $g_i\in \mathbb{F}_q$ generates a
cyclic subgroup of $\mathbb{F}_q\sim H_1$
which is invariant under the
action of the group $\mathbb{A}_{q,m}$. The flag of $n/r$
linear subspaces of $\mathbb{F}_q$ determined by this
basis defines a nested sequence
$H_{n/r}\triangleleft\dots\triangleleft H_1$ of normal subgroups
of $H_1$, and  the set $g_1,g_1^{-1},\dots,
g_{n/q},g_{n/q}^{-1},\xi^m,\xi^{-m}$ is a
$(G_{q,m},H_1,\dots,H_m)$-admissible set of generators for
$\Gamma$.

Now since $H_1$ equals the commutator of $G_{q,m}$, every
character of $G_{q,m}$ (i.e. a one-dimensional unitary
representation of $G_{q,m}$) factors to a character of
$G_{q,m}/H_1$. On the other hand, it is well known [M] that the
dimension of every irreducible representation of $G_{q,m}$ which
is not a character is at least $(q-1)/m$. Thus by Corollary 4.5
the group $G_{q,m}$ gives rise to manifolds for which the $j$-th
eigenvalue for $2\leq j\leq n/r$ has multiplicity at least
$(q-1)/m$. Since $n/r$ and $(q-1)/m$ can be chosen arbitrarily
large our Theorem $C$ from the introduction follows.

\bigskip

\noindent{\bf Acknowledgement;} We like to thank the anonymous
referee for pointing out a mistake in an earlier version of this
paper. The first author is also grateful to Dorothee Schueth for
helpful discussions.

\section{References}

\begin{enumerate}
\item[{[BC]}] M. Burger, B. Colbois,
{\it A propos de la multiplicit\'e de la premi\`ere valeur propre
d'une surface de Riemann}, C.R.A.S. Paris 300 (1985), 247-249.
\item[{[B]}] P. Buser, {\sl Geometry and spectra of compact Riemann surfaces},
\hfil\break Birkh\"auser 1992.
\item[{[BCD]}] P. Buser, B. Colbois, J. Dodziuk,
{\it Tubes and eigenvalues for
negatively curved manifolds}, J. Geom. Anal. 3 (1993), 1-26.
\item[{[CC1]}] B. Colbois, G. Courtois, {\it Les valeurs
propers inf\'erieures
\`a $\frac{1}{4}$ des surfaces de riemann de petit rayon
d'injectivit\'e}, Comm. Math. Helv. 64 (1989), 349-362.
\item[{[CC2]}] B. Colbois, G. Courtois,
{\it Convergence de vari\'et\'es et convergence du spectre du
laplacien}, Ann. Sc. Ec. Norm. Sup. 4 (1991), 507-518.
\item[{[D]}] E.~B.~Davies, {\sl Spectral theory and differential
operators}, Cambridge University Press, Cambridge 1995.
\item[{[F]}] K. Fissmer, {\it Diskretes Spektrum von hyperbolischen Fl\"achen:
Beispiele f\"ur Konvergenz}, Diplomarbeit, Bonn 2000.
\item[{[G]}] M. Gromov, {\sl Metric structures for Riemannian
and non-Riemannian spaces}, Birkh\"auser 1999.
\item[{[dH]}] P. de la Harpe, {\sl Topics in geometric group theory},
Univ. of Chicago Press, Chicago 2000.
\item[{[J]}] C. Judge, {\it Tracking eigenvalues to the frontier of moduli
space I: Convergence and spectral accumulation}, J. Funct. Analysis 184 (2001), 273-290.
\item[{[M]}] G.W. Mackey, {\sl Induced representations of groups and
quantum mechanics}, Benjamin Inc. 1968.

\end{enumerate}

\end{document}